\documentclass[12pt]{article}
\usepackage{amsfonts}
\usepackage{latexsym,amsmath,amssymb}
\usepackage{graphics,epstopdf}
\usepackage[pdftex]{graphicx}

\newtheorem{definition}{Definition}[section]
\newtheorem{theorem}[definition]{Theorem}

\newtheorem{lemma}[definition]{Lemma}
\numberwithin{equation}{section}
\begin{document}

\begin{center}
{\Large \textbf{\ Approximation by Kantorovich type $(p,q)$-Bernstein-Schurer Operators }}

\bigskip

\textbf{M. Mursaleen}, \textbf{Faisal Khan}

Department of\ Mathematics, Aligarh Muslim University, Aligarh--202002, India%
\\[0pt]

mursaleenm@gmail.com; faisalamu2011@gmail.com \\[0pt]

\bigskip

\bigskip

\textbf{Abstract}
\end{center}

\parindent=8mm {\footnotesize {In this paper, we introduce a Shurer type genaralization of
$(p,q)$-Bernstein-Kantorovich operators based on $(p,q)$-integers and we call it as
$(p,q)$-Bernstein-Schurer Kantorovich operators. We study approximation properties for
these operators based on Korovkin's type approximation theorem and also study some direct
theorems. Furthermore, we give comparisons and some illustrative graphics for the convergence
of operators to some function.}}

\bigskip



\section{Introduction and preliminaries}

\hspace{8mm}The applications of $q$-calculus emerged as a new area in the field of approximation
theory from last two decades. The development of $q$-calculus has led to the discovery of
various modifications of Bernstein polynomials involving $q$-integers. The aim of these
generalizations is to provide appropriate and powerful tools to application
areas such as numerical analysis, computer-aided geometric design and
solutions of differential equations.

\parindent8mmIn 1987, Lupa\c{s} \cite{lp} introduced the first $q$-analogue of Bernstein
operators \cite{brn} and investigated its approximating and shape-preserving properties.
Another $q$-generalization of the classical Bernstein polynomials is due to Phillips \cite{pl}.
Several generalizations of well-known positive linear operators based on $q$-integers were
introduced and their approximation properties have been studied by several authors. For instance,
the approximation properties of $q$-Bleimann, Butzer and Hahn operators \cite{ar1}; $q$-analogue
of Sz\'{a}sz-Kantorovich operators \cite{m3}; Approximation by Kantorovich type $q$-Bernstein
operators \cite{dal} and $q$-analogue of generalized Berstein-Shurer operators {\cite{ma1} were studied.

Mursaleen et al studied approximation properties of the $q$-analogue of generalized Berstein
Shurer operators  {\cite{ma1}.\\

 Recently, Mursaleen et al applied $(p,q)$-calculus in approximation theory and introduced $(p,q)$-analogue of Bernstein Operators \cite{mka1} and $(p,q)$-analogue of Bernstein-Stancu Operators \cite{mka2}, $(p,q)$-analogue of
Bernstein-Kantorovich Operators \cite{mka3} respectively. Recently these works got attention in computer aided geometric design (CAGD)\\

   As in computer aided geometric design , the basis of Bernstein polynomials plays a significant role in order to preserve the shape of the curves or surfaces. The classical Bezier curve \cite{bezier} and $q$-bezier \cite{ph1,hp} constructed with Bernstein basis functions are the most important curve in CAGD. Thus motivated by the work of Mursaleen et al in \cite{mka1}, very recently, Khalid et al in \cite{kl} introduced  Bezier curves and surfaces in Computer aided geometric design defined by $(p,q)$-integers.\\

Ozarslan and Vedi \cite{ov} introduced  Kantorovich type genelaization of the $q$-Bernstein-Schurer
operators and studied some approximation properties of these operators. These operators are given as follows:
\begin{equation}\label{e1.1}
K^q_{n,\ell}(f;x)=\sum\limits_{k=0}^{n+\ell}b^q_{n,\ell,k}(x)~\int_{0}^{1}f\left(\frac{[k]_{q}}{[n+1]_{q}}+
\frac{1+(q-1)[k]_{q}}{[n+1]_q}t\right)d_{q}t ,~~x\in \lbrack 0,1],
\end{equation}
\begin{equation*}
  b^q_{n,\ell,k}(x):= \left[
\begin{array}{c}
n+\ell \\
k
\end{array}
\right] _{q}x^{k}\prod\limits_{s=0}^{n+\ell-k-1}(1-q^{s}x).
\end{equation*}
Where $\ell\in\mathbb{N_{\circ}}=\mathbb{N}\cup\{0\},~~0<q<1$ is fixed and
$K_{n,\ell}^q:C[0,1+\ell]\rightarrow C[0,1+\ell]$ are defined for all $n\in\mathbb{N}$
and for any function $f\in C[0,1+\ell].$\\

\parindent=8mm Details on the $q$-calculus can be found in \cite{kac} and for the applications of
$q$-calculus in approximation theory, one can refer \cite{gp1}.\\

Mursaleen et al introduced $(p,q)$-Bernstein-kantorovich operators as follows in \cite{mka3}

\begin{equation*}
K_{n}^{(p,q)}(f;x)=[n+1]_{p,q}\sum\limits_{k=0}^{n}\frac{(p-q)}{p^k(p-1)-q^k(k-1)} b_{n,k}^{(p,q)}(x)\int_{\frac{[k]_{p,q}}{[n+1]_{p,q}}}^{\frac{[k+1]_{p,q}}{[n+1]_{p,q}}}f(t)d_{p,q}t ,~~x\in \lbrack 0,1]
\end{equation*}
where
\begin{equation*}
   b_{n,k}^{(p,q)}(x)= \left[
\begin{array}{c}
n \\
k%
\end{array}%
\right] _{p,q}x^k(1-x)_{p,q}^{n-k}= \left[
\begin{array}{c}
n \\
k%
\end{array}%
\right] _{p,q}x^{k}\prod\limits_{s=0}^{n-k-1}(p^s-q^{s}x).
\end{equation*}

Motivated by the work of  Mursaleen et al
for $(p,q)$-Bernstein-Kantorovich  operators based on
$(p,q)$-integers in  \cite{mka3}, now we present a Shurer type generalisation of $(p,q)$-Bernstein-Kantorovich  operators.\\

Let us recall certain notations of $(p,q)$-calculus:

The $(p,q)$-integer was introduced in order to generalize or unify several
forms of $q$-oscillator algebras well known in the earlier physics
literature related to the representation theory of single parameter quantum
algebras \cite{chak}. The $(p,q)$-integer $[n]_{p,q}$ is defined by

\begin{equation*}
[n]_{p,q}:=\frac{p^n-q^{n}}{p-q},~~~n=0,1,2,\cdots,~~0<q<p\leq1.
\end{equation*}

The $(p,q)$-Binomial expansion is
\begin{equation*}
(ax+by)_{p,q}^n :=\sum\limits_{k=0}^{n}\left[
\begin{array}{c}
n \\
k%
\end{array}%
\right] _{p,q}q^{\frac{k(k-1)}{2}}p^{\frac{(n-k)(n-k-1)}{2}}a^{n-k}b^kx^{n-k}y^k
\end{equation*}

\begin{equation*}
(x+y)_{p,q}^n:=(x+y)(px+qy)(p^2x+q^2y)\cdots(p^{n-1}x+q^{n-1}y).
\end{equation*}

Also, the $(p,q)$-binomial coefficients are defined by

\begin{equation*}
\left[
\begin{array}{c}
n\\
k%
\end{array}%
\right] _{p,q}:=\frac{[n]_{p,q}!}{[k]_{p,q}![n-k]_{p,q}!}
\end{equation*}%
and the definite integrals of the function $f$ are defined by
\begin{equation*}
    \int_{0}^{a}f(x)d_{p,q}x=(q-p)a\sum\limits_{k=0}^{\infty}\frac{p^k}{q^{k+1}}
    f\left(\frac{p^k}{q^{k+1}}a\right),~~~\text{when}~~~\left|\frac pq\right|<1,
\end{equation*}
and
\begin{equation*}
    \int_{0}^{a}f(x)d_{p,q}x=(p-q)a\sum\limits_{k=0}^{\infty}\frac{q^k}{p^{k+1}}
    f\left(\frac{q^k}{p^{k+1}}a\right),~~~\text{when}~~~\left|\frac pq\right|>1.
\end{equation*}

Details on $(p,q)$-calculus can be found in \cite{mah,jacob,jag,vivek,sad}.
For $p=1$, all the notions of $(p,q)$-calculus are reduced to $q$-calculus.

\parindent=8mm Now, we introduce a new generalization of $q$-Bernstein-Kantorovich
operators by using the notion of $(p,q)$-calculus and call it as $(p,q)$-Bernstein-Schurer-Kantorovich  operators.
We study the approximation properties based on Korovkin's type approximation theorem and also establish
some direct theorems. Further, we show comparisons and some illustrative graphics for the convergence of
operators to a function.\\

\section{Construction of Operators}

Now, we introduce $(p,q)$-analogue of Kantorovich type Bernstein-Schurer operators as
\begin{equation}\label{e2.1}
K_{n,\ell}^{(p,q)}(f;x)=\sum\limits_{k=0}^{n+\ell} b_{n,\ell,k}^{(p,q)}(x)~\int_{0}^{1}f\left(\frac{[k]_{p,q}}{[n+1]_{p,q}}+
\frac{[k+1]_{p,q}-[k]_{p,q}}{[n+1]_q}t\right)d_{p,q}t ,~~x\in \lbrack 0,1],
\end{equation}
where
\begin{equation*}
   b_{n,\ell,k}^{(p,q)}(x)= \left[
\begin{array}{c}
n+\ell \\
k%
\end{array}
\right] _{p,q}x^k(1-x)_{p,q}^{n+\ell-k}= \left[
\begin{array}{c}
n+\ell \\
k%
\end{array}
\right] _{p,q}x^{k}\prod\limits_{s=0}^{n+\ell-k-1}(p^s-q^{s}x).
\end{equation*}

For $p=1$, equation (\ref{e2.1}) turns out to be the classical $q$-Bernstein-Kantorovich operators (\ref{e1.1}).\\
\\
Now, we have the following basic lemmas:
 \begin{lemma}\label{lem2.1} For $x\in \lbrack 0,1],~0<q<p\leq 1$
\begin{enumerate}
\item [(i)] $K_{n,\ell}^{(p,q)}(1;x)=~1$;
\item[(ii)]$K_{n,\ell}^{(p,q)}(t;x)=~\frac{(px+1-x)_{p,q}^{n+\ell}}{[2]_{p,q}[n+1]_{p,q}}
+\frac{(p+2q-1)[n+\ell]_{p,q}}{[2]_{p,q}[n+1]_{p,q}}x$;
 \item[(iii)] $K_{n,\ell}^{(p,q)}(t^2;x)=~\frac{(p^2x+1-x)_{p,q}^{n+\ell}}{[3]_{p,q}[n+1]_{p,q}^2}
+\left(1+\frac{2q}{[2]_{p,q}}+\frac{q^2-1}{[3]_{p,q}}\right)\frac{[n+\ell]_{p,q}}{[n+1]_{p,q}^2}(px+1-x)_{p,q}^{n+\ell-1}x\\
~~~~~~~~~~~~~~~~~~~~+\left(1+\frac{2(q-1)}{[2]_{p,q}}+\frac{(q-1)^2}{[3]_{p,q}}\right)\frac{[n+\ell]_{p,q}[n+\ell-1]_{p,q}}{[n+1]_{p,q}^2}x^2;$\\
 \item[(iv)]$K_{n,\ell}^{(p,q)}\big{(}(t-x);x\big{)}=\frac{(p^2x+1-x)_{p,q}^{n+\ell}}{[2]_{p,q}[n+1]_{p,q}}+\left(\frac{(p+2q-1)}{[2]_{p,q}[n+1]_{p,q}}-1\right)x;$
\item[(v)]$K_{n,\ell}^{(p,q)}\big{(}(t-x)^2;x\big{)}=\frac{(p^2x+1-x)_{p,q}^{n+\ell}}{[3]_{p,q}[n+1]_{p,q}^2}
~+~\left\{\left(1+\frac{2q}{[2]_{p,q}}+\frac{q^2-1}{[3]_{p,q}}\right)\frac{[n+\ell]_{p,q}(px+1-x)_{p,q}^{n+\ell-1}}{[n+1]_{p,q}^2}
-\frac{2(px+1-x)_{p,q}^{n+\ell}}{[2]_{p,q}[n+1]_{p,q}}\right\}x\\
~~~~~~~~~~~~~~~~~~~~~~~~~~~~~~~~+\left\{q\left(1+\frac{2(q-1)}{[2]_{p,q}}+\frac{(q-1)^2}{[3]_{p,q}}\right)\frac{[n+\ell]_{p,q}[n+\ell-1]_{p,q}}{[n+1]_{p,q}^2}
-\frac{2(p+2q-1)[n+\ell]_{p,q}}{[2]_{p,q}[n+1]_{p,q}}+1\right\}x^2$.
\end{enumerate}
\end{lemma}

\section{Approximation results}

\hspace{8mm}In this section, we find out the rate of convergence of the operators (\ref{e2.1})
using the modulus of continuity and Lipschitz classes. Furthermore, we clculate the rate convergence
in terms of the first modulus of continuity of the function.

\parindent=8mmLet $C[0,1+\ell]$ be the linear space of all real valued continuous functions $f$
on $[0,1+\ell]$ and let $T$ be a linear operator defined on $C[0,1+\ell].$
We say that $T$ is $positive$ if for every non-negative $f\in C[0,1+\ell],$ we
have $T(f,x)\geq 0$ for all $x\in $ $[0,1+\ell].$

\parindent=8mmThe classical Korovkin approximation theorem \cite{alt, pp,sr} states as follows:

\parindent=8mmLet $(T_{n})$ be a sequence of positive linear operators from
$C[0,1+\ell]$ into $C[a,b].$ Then $\lim_{n}\Vert T_{n}(f,x)-f(x)\Vert_{C[0,1+\ell]}=0$,
for all $f\in C[0,1+\ell]$ if and only if $\lim_{n}\Vert T_{n}(f_{i},x)-f_{i}(x)\Vert _{C[0,1+\ell]}=0$,
for $i=0,1,2$, where $f_{0}(x)=1,~f_{1}(x)=x$ and $f_{2}(x)=x^{2}.$

\begin{theorem}\label{theorem3.1} {\em Let $0<q_{n}<p_{n}\leq 1$ such that
$\lim\limits_{n\rightarrow \infty }p_{n}=1$ and $\lim\limits_{n\rightarrow
\infty }q_{n}=1$. Then for each $f\in [0,1+\ell],~K_{n,\ell}^{(p_{n},q_{n})}(f;x)$
converges uniformly to $f$ on $[0,1+\ell]$.}
\end{theorem}

\parindent=0mm\textbf{Proof}. By the Korovkin Theorem it is sufficient to
show that
\begin{equation*}
\lim\limits_{n\rightarrow \infty }\Vert
K_{n,\ell}^{(p_{n},q_{n})}(t^{m};x)-x^{m}\Vert _{C[0,1+\ell]}=0,~~~m=0,1,2.
\end{equation*}%
By Lemma \ref{lem2.1} (i), it is clear that
\begin{equation*}
\lim\limits_{n\rightarrow \infty }\Vert K_{n,\ell}^{(p_{n},q_{n})}(1;x)-1\Vert
_{C[0,1+\ell]}=0.
\end{equation*}%
Now, by Lemma \ref{lem2.1}  (ii)
\begin{equation*}
 |K_{n,\ell}^{(p_{n},q_{n})}(t;x)-x|\leq \frac{(p_nx+1-x)_{p_n,q_n}^{n+\ell}}{[2]_{p,q}[n+1]_{p_n,q_n}}
 +\left(\frac{(p_n+2q_n-1)}{[2]_{p_n,q_n}[n+1]_{p_n,q_n}}-1\right)x
\end{equation*}
Taking maximum of both sides of the above inequality, we get
\begin{equation*}
    \|K_{n,\ell}^{(p_{n},q_{n})}(t;x)-x\|_{[0,1+\ell]}\leq\frac{p_n^n}{[2]_{p_n,q_n}[n+1]_{p_n,q_n}}
    +\frac{(p_n+2q_n-1)[n+\ell]_{p_n,q_n}}{[2]_{p_n,q_n}[n+1]_{p_n,q_n}}-1
\end{equation*}
which yields
\begin{equation*}
\lim\limits_{n\rightarrow \infty }\Vert
K_{n,\ell}^{(p_{n},q_{n})}(t;x)-x\Vert _{C[0,1+\ell]}=0.
\end{equation*}%
Similarly we can show that
\begin{equation*}
\lim\limits_{n\rightarrow \infty }\Vert
K_{n,\ell}^{(p_{n},q_{n})}(t^2;x)-x^2\Vert _{C[0,1+\ell]}=0.
\end{equation*}%
Thus the proof is completed.\newline
\newline

Now we will compute the rate of convergence in terms of modulus of
continuity.\newline

\parindent=8mm Let $f\in [0,1+\ell]$. The modulus of continuity of $f$ denoted
by $\omega_{f}(\delta )$ gives the maximum oscillation of $f$ in any interval of
length not exceeding $\delta >0$ and it is given by the relation
\begin{equation*}
\omega_{f}(\delta )=\sup\limits_{|x-y|\leq \delta }|f(x)-f(y)|,~~x,y\in [0,1+\ell].
\end{equation*}%
It is known that $\lim\limits_{\delta \rightarrow 0^+}\omega_{f}(\delta )=0$ for $%
f\in C[0,1+\ell]$ and for any $\delta >0$ one has
\begin{equation}\label{e3.1}
|f(y)-f(x)|\leq \omega_{f}(\delta)\biggl{(}\frac{|y-x|}{\delta}+1%
\biggl{)}.
\end{equation}

\begin{theorem}\label{theorem3.2} {\em  If $f\in C[0,1]$, then
\begin{equation*}
\bigl{|}K_{n,\ell}^{(p,q)}(f;x)-f(x)\bigl{|}\leq 2\omega\left(f,\sqrt{\delta _{n,\ell}^{(p,q)}(x)}\right)
\end{equation*}%
takes place, where $\delta_{n,\ell}^{(p,q)}=K_{n,\ell}^{(p,q)}\big{(}(t-x)^2;x\big{)}$}
\end{theorem}

\parindent=0mm \textbf{Proof}. Since $K_{n,\ell}^{(p,q)}(1;x)=1$, we have

\begin{align*}
\big|K_{n,\ell}^{(p,q)}(f;x)-f(x)\big|&\leq K_{n,\ell}^{(p,q)}\big(|f(t)-f(x)|;x\big)\\
&=\sum\limits_{k=0}^{n+\ell} b_{n,\ell,k}^{(p,q)}(x)~\int_{0}^{1}f\left|\left(\frac{[k]_{p,q}}{[n+1]_{p,q}}+
\frac{[k+1]_{p,q}-[k]_{p,q}}{[n+1]_q}t\right)-f(x)\right|d_{p,q}t.
\end{align*}
In view of (\ref{e3.1}), we get
\begin{align*}
  \big|K_{n,\ell}^{(p,q)}(f;x)-f(x)\big|&\leq\sum\limits_{k=0}^{n+\ell} b_{n,\ell,k}^{(p,q)}(x)~\int_{0}^{1}
\left(\frac{\left|\frac{[k]_{p,q}}{[n+1]_{p,q}}+\frac{[k+1]_{p,q}-[k]_{p,q}}{[n+1]_q}t-x\right|}{\delta}+1\right)\omega(f,\delta)~d_{p,q}t\\
&=\omega(f,\delta)+\frac{\omega(f,\delta)}{\delta}\sum\limits_{k=0}^{n+\ell} b_{n,\ell,k}^{(p,q)}(x)
\int_{0}^{1}\left|\frac{[k]_{p,q}}{[n+1]_{p,q}}+\frac{[k+1]_{p,q}-[k]_{p,q}}{[n+1]_q}t-x\right|d_{p,q}t
\end{align*}
Now using Cauchy-Schwartz inequality, we get
\begin{align*}
\big|K_{n,\ell}^{(p,q)}(f;x)-f(x)\big|&\leq\omega(f,\delta)+\frac{\omega(f,\delta)}{\delta}\sum\limits_{k=0}^{n+\ell} b_{n,\ell,k}^{(p,q)}(x)\\
&\hspace{3cm}\times\bigg\{\int_{0}^{1}\bigg(\frac{[k]_{p,q}}{[n+1]_{p,q}}+\frac{[k+1]_{p,q}-[k]_{p,q}}{[n+1]_q}t-x\bigg)^2d_{p,q}t\bigg\}^{\frac{1}{2}}.
\end{align*}

Again applying the cauchy-Schwartz inequality, we have

\begin{equation*}
\big|K_{n,\ell}^{(p,q)}(f;x)-f(x)\big|\leq\omega(f,\delta)+\frac{\omega(f,\delta)}{\delta}\bigg\{\sum\limits_{k=0}^{n+\ell} b_{n,\ell,k}^{(p,q)}(x)\hspace{5.5cm}
\end{equation*}
\begin{equation*}
\times\int_{0}^{1}\bigg(\frac{[k]_{p,q}}{[n+1]_{p,q}}+\frac{[k+1]_{p,q}-[k]_{p,q}}{[n+1]_q}t-x\bigg)^2d_{p,q}t\bigg\}^{\frac{1}{2}}
\bigg\{\sum\limits_{k=0}^{n+\ell} b_{n,\ell,k}^{(p,q)}(x)\bigg\}^{\frac{1}{2}}
\end{equation*}

\begin{align*}
\hspace{1.5cm}&=\omega(f,\delta)+\frac{\omega(f,\delta)}{\delta}\bigg\{\sum\limits_{k=0}^{n+\ell}b_{n,\ell,k}^{(p,q)}(x)
\int_{0}^{1}\bigg(\frac{[k]_{p,q}}{[n+1]_{p,q}}+\frac{[k+1]_{p,q}-[k]_{p,q}}{[n+1]_q}t-x\bigg)^2d_{p,q}t\bigg\}^{\frac{1}{2}}
\end{align*}
\begin{equation*}
=\omega(f,\delta)+\frac{\omega(f,\delta)}{\delta}\bigg\{K_{n,\ell}^{(p,q)}\big{(}(t-x)^2;x\big{)}\bigg\}^{\frac{1}{2}}.\hspace{4cm}
\end{equation*}
Now on taking $\delta:=\delta_{n,\ell}^{(p,q)}(x)=K_{n,\ell}^{(p,q)}\big{(}(t-x)^2;x\big{)},$ we obtain
\begin{equation*}
  \big|K_{n,\ell}^{(p,q)}(f;x)-f(x)\big|\leq2\omega\left(f,\sqrt{K_{n,\ell}^{(p,q)}\big((t-x)^2;x\big)}\right)
\end{equation*}
 This completes the proof of the theorem.\\

\hspace{8mm}Now we give the rate of convergence of the operators $K_{n,\ell}^{(p,q)}$ in terms
of the elements of the usual Lipschitz class $\text{Lip}_M(\alpha)$.\\

\hspace{8mm}Let $f\in C[0,1+\ell]$, $M>0$ and $0<\alpha\leq1$.We recall that $f$ belongs to the class
$\text{Lip}_M(\alpha)$ if the inequality
\begin{equation*}
    |f(t)-f(x)|\leq M|t-x|^\alpha~~~(t,x\in[0,1+\ell])
\end{equation*}
is satisfied.

\begin{theorem}\label{theorem3.3} {\em Let $0<q<p\leq1$. Ten for each $f\in \text{Lip}_M(\alpha)$ we have
\begin{equation*}
    |K_{n,\ell}^{(p,q)}(f;x)-f(x)|\leq M\left(\delta_{n,\ell}^{(p,q)}(x)\right)^{\frac{\alpha}{2}}
\end{equation*}
where
\begin{equation*}
\delta_{n,\ell}^{(p,q)}(x)={K_{n,\ell}^{(p,q)}\big{(}(t-x)^2;x\big{)}}.
\end{equation*}}
\end{theorem}

\parindent=0mm\textbf{Proof}. By the monotonicity of the operators $K_{n,\ell}^{(p,q)}$, we can write
\begin{align*}
\big|K_{n,\ell}^{(p,q)}(f;x)-f(x)\big|&\leq K_{n,\ell}^{(p,q)}\big(|f(t)-f(x)|;x\big)\\
&=\sum\limits_{k=0}^{n+\ell} b_{n,\ell,k}^{(p,q)}(x)~\int_{0}^{1}f\left|\left(\frac{[k]_{p,q}}{[n+1]_{p,q}}+
\frac{[k+1]_{p,q}-[k]_{p,q}}{[n+1]_q}t\right)-f(x)\right|d_{p,q}t\\
 &\leq M\sum\limits_{k=0}^{n+\ell} b_{n,\ell,k}^{(p,q)}(x)~\int_{0}^{1}\left|\left(\frac{[k]_{p,q}}{[n+1]_{p,q}}+
\frac{[k+1]_{p,q}-[k]_{p,q}}{[n+1]_q}t\right)-x\right|^{\alpha}d_{p,q}t\\.
\end{align*}
Now applying the H\"{o}lder's inequality for the sum with $p=\frac2{\alpha}$ and $q=\frac2{2-\alpha}$ and taking into consideration Lemma 2.1(i) and Lemma 2.2(ii), we have
\begin{eqnarray*}
  |K_{n,\ell}^{(p,q)}(f;x)-f(x)| &\leq& M\sum\limits_{k=0}^{n+\ell}\Bigg{\{} b_{n,\ell,k}^{(p,q)}(x)\int_{0}^{1}\left(\left(\frac{[k]_{p,q}}{[n+1]_{p,q}}+
\frac{[k+1]_{p,q}-[k]_{p,q}}{[n+1]_q}t\right)-x\right)^{2}d_{p,q}t \Bigg{\}}^{\frac\alpha2}\\
  &&\times\Bigg{\{} b_{n,\ell,k}^{(p,q)}(x)\int_{0}^{1} 1 d_{p,q}t \Bigg{\}}^{\frac{2-\alpha}2}\\
 &\leq& M\Bigg{\{}\sum\limits_{k=0}^{n+\ell} b_{n,\ell,k}^{(p,q)}(x)\int_{0}^{1}\left(\left(\frac{[k]_{p,q}}{[n+1]_{p,q}}+
\frac{[k+1]_{p,q}-[k]_{p,q}}{[n+1]_q}t\right)-x\right)^{2}d_{p,q}t \Bigg{\}}^{\frac\alpha2}\\
  &&\times\Bigg{\{}\sum\limits_{k=0}^{n+\ell} b_{n,\ell,k}^{(p,q)}(x)\int_{0}^{1} 1 d_{p,q}t \Bigg{\}}^{\frac{2-\alpha}2}\\
  &=&M\bigl{\{}K_{n,\ell}^{(p,q)}\bigl{(}(t-x)^2;x\bigl{)}\bigl{\}}^\frac{\alpha}{2}.
\end{eqnarray*}
Choosing $\delta_{n,\ell}^{(p,q)}(x)={K_{n,\ell}^{(p,q)}\big{(}(t-x)^2;x\big{)}}$, we arrive at the desired result.
\\

\parindent=8mm Next we prove the local approximation property for the operators $K_{n,\ell}^{(p,q)}$.
The Peetre's $K$-functional is defined by \newline
\begin{equation*}
K_{2}(f,\delta )=\inf [\{\Vert f-g\Vert +\delta \Vert g^{\prime \prime}\Vert \}:g\in W^{2}],
\end{equation*}%
where%
\begin{equation*}
W^{2}=\{g\in C[0,1+\ell]:g^{\prime },g^{\prime \prime }\in C[0,1+\ell]\}.
\end{equation*}

By \cite{dl}, there exists a positive constant $C>0$ such that $%
K_{2}(f,\delta )\leq C\omega_{2}(f,\delta ^{\frac{1}{2}}),~~\delta >0;$ where the
second order modulus of continuity is given by
\begin{equation*}
\omega_{2}(f,\delta ^{\frac{1}{2}})=\sup\limits_{0<h\leq \delta ^{\frac{1}{2}%
}}\sup\limits_{x\in \lbrack 0,1+\ell]}\mid f(x+2h)-2f(x+h)+f(x)\mid .
\end{equation*}%
\newline

Also for $f\in [0,1+\ell]$ the usual modulus of continuity is given by
\begin{equation*}
\omega(f,\delta )=\sup\limits_{0<h\leq \delta ^{\frac{1}{2}}}\sup\limits_{x\in
\lbrack 0,1+\ell]}\mid f(x+h)-f(x)\mid .
\end{equation*}

\begin{theorem}\label{theorem3.4} {\em Let $f\in C[0,1+\ell]$ and $0<q<p\leq1$.
Then for all $n\in \mathbb{N}$, there exists an absolute constant $C>0$
such that
\begin{equation*}
\mid K_{n,\ell}^{(p,q)}(f;x)-f(x)\mid \leq C\omega_{2}\big(f,\sqrt{a_{n,\ell}^{(p,q)}(x)}\big)
+\omega\big(f,c_{n,\ell}^{(p,q)}\big),
\end{equation*}%
where%
\begin{equation*}
a_{n,\ell}^{(p,q)}(x)= K_{n,\ell}^{(p,q)}\left( (t-x)^{2};x\right) +\left(\alpha_{n,\ell}^{(p,q)}(x)-x\right)^2,~{c_{n,\ell}^{(p,q)}(x)=\left(\alpha_{n,\ell}^{(p,q)}(x)-x\right)}
\end{equation*}
and
\begin{equation*}
 \alpha_{n,\ell}^{(p,q)}(x)=\frac{(px+1-x)_{p,q}^{n+\ell}}{[2]_{p,q}[n+1]_{p,q}}
+\frac{(p+2q-1)[n+\ell]_{p,q}}{[2]_{p,q}[n+1]_{p,q}}x.
\end{equation*}}
\end{theorem}

\parindent=0mm \textbf{Proof}. Using the operators (\ref{e2.1}), we define the following operators
\begin{equation}\label{e3.2}
  \tilde{K}_{n,\ell}^{(p,q)}(f;x):=K_{n,\ell}^{(p,q)}(f;x)-f\bigg(\frac{(px+1-x)_{p,q}^{n+\ell}}{[2]_{p,q}[n+1]_{p,q}}
+\frac{(p+2q-1)[n+\ell]_{p,q}}{[2]_{p,q}[n+1]_{p,q}}x\bigg)+f(x).
\end{equation}

 Then, by the Lemma \ref{lem2.1} (ii) For $x\in \lbrack 0,1],~0<q<p\leq 1$
\begin{enumerate}
\item [] $\tilde{K}_{n,\ell}^{(p,q)}(1;x)=~1$;
\item[]$\tilde{K}_{n,\ell}^{(p,q)}(t-x;x)=0$.
\end{enumerate}
Then, for a given $g\in W^{2}.$ From Taylor's expansion, we get
\begin{equation*}
g(t)=g(x)+g^{\prime }(x)(t-x)+\int_{x}^{t}(t-u)~g^{\prime \prime
}(u)~du,~~~t\in \lbrack 0,A],~A>0,~~
\end{equation*}
and for the operators (\ref{e3.2}), we get

\begin{eqnarray*}
\tilde{K}_{n,\ell}^{(p,q)}(g;x)&=&g(x)+\tilde{K}_{n,\ell}^{(p,q)}\left( \int_{x}^{t}(t-u)~g^{\prime \prime
}(u)~du;x\right),\\
 \mid \tilde{K}_{n,\ell}^{(p,q)}(g;x)-g(x)\mid& =& \biggl{|}\tilde{K}_{n,\ell}^{(p,q)}\left(
\int_{x}^{t}(t-u)~g^{\prime \prime }(u)~du;x\right) \biggl{|}\\
{}&=&\left|{K}_{n,\ell}^{(p,q)}\left(\int_{x}^{t}(t-u)~g^{\prime \prime }(u)~du;x\right)-\int_{x}^{\alpha_{n,\ell}^{(p,q)}(x)}\left(\alpha_{n,\ell}^{(p,q)}(x)-u\right)g^{\prime\prime}(u)~du\right|\\
{}&\leq&\left|{K}_{n,\ell}^{(p,q)}\left(\int_{x}^{t}(t-u)~g^{\prime \prime }(u)~du;x\right)\right|+\left|\int_{x}^{\alpha_{n,\ell}^{(p,q)}(x)}\left(\alpha_{n,\ell}^{(p,q)}(x)-u\right)g^{\prime\prime}(u)~du\right|\\
&\leq& \Vert g^{\prime \prime }\Vert K_{n,\ell}^{(p,q)}\left( (t-x)^{2};x\right) +\|g^{\prime\prime}\|\left(\alpha_{n,\ell}^{(p,q)}(x)-x\right)^2\\
&=&\Vert g^{\prime \prime }\Vert\left( K_{n,\ell}^{(p,q)}\left( (t-x)^{2};x\right) +\left(\alpha_{n,\ell}^{(p,q)}(x)-x\right)^2\right)
= \Vert g^{\prime \prime}\Vert a_{n,\ell}^{(p,q)}(x).
\end{eqnarray*}

On the other hand, by the definition of $\tilde{K}_{n,\ell}^{(p,q)}(f;x),$ we have
\begin{equation*}
\mid \tilde{K}_{n,\ell}^{(p,q)}(f;x)\mid \leq 4\Vert f\Vert .
\end{equation*}
Now
\begin{eqnarray*}
\mid K_{n,\ell}^{(p,q)}(f;x)-f(x)\mid &\leq& \mid \tilde{K}_{n,\ell}^{(p,q)}\bigl{(}(f-g);x\bigl{)}
-(f-g)(x)\mid +\mid \tilde{K}_{n,\ell}^{(p,q)}(g;x)-g(x)\mid\\
  &&+\left|f\big(\alpha_{n,\ell}^{(p,q)}(x)\big)-f(x)\right|\\
&\leq& 4\Vert f-g\Vert +a_{n,\ell}^{(p,q)}(x)\Vert g^{\prime \prime }\Vert+\left|f\big(\alpha_{n,\ell}^{(p,q)}(x)\big)-f(x)\right|\\
&\leq&4(\Vert f-g\Vert +a_{n,\ell}^{(p,q)}(x)\Vert g^{\prime \prime} \Vert)+\omega\big(f,c_{n,\ell}^{(p,q)}\big).
\end{eqnarray*}

In view of the property of $K$-functional, it can be easily seen that

\begin{eqnarray*}
\mid K_{n,\ell}^{(p,q)}(f;x)-f(x)\mid &\leq& 4 K\big(f,\sqrt{a_{n,\ell}^{(p,q)}(x)}\big)
+\omega\big(f,c_{n,\ell}^{(p,q)}\big)\\
&\leq & C\omega_{2}\big(f,\sqrt{a_{n,\ell}^{(p,q)}(x)}\big)
+\omega\big(f,c_{n,\ell}^{(p,q)}\big).
\end{eqnarray*}

\parindent=8mm This completes the proof of the theorem.


\newpage

\section{Graphical analysis}

With the help of Matlab, we show comparisons and some illustrative graphics
\cite{ma1} for the convergence of $(p,q) $-Bernstein-Schurer-Kantorovich operators
$K_{n,\ell}^{(p,q)}$ to the function $f (x) = 1+cos(5x^2),$ for different values of
parameters $p,q,n$, convergence of the operators $K_{n,\ell}^{(p,q)}$ to the function
is shown in figure $1,2,3,~\text{and}~4.$

\begin{figure*}[htb!]
\begin{center}
\includegraphics[height=6cm, width=10cm]{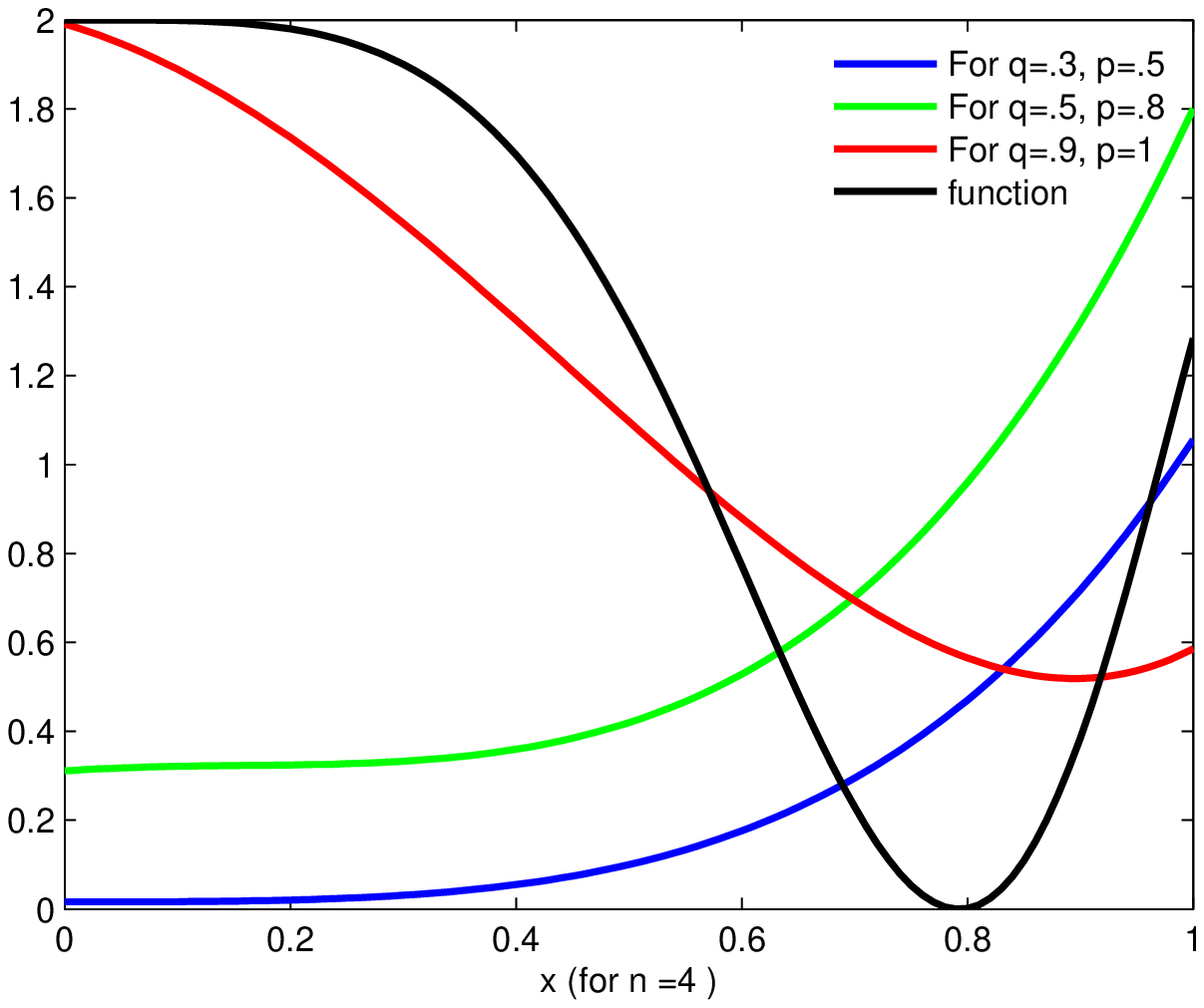}
\end{center}
\caption{ Approximation by $(p,q)$-Bernstein-Schurer-Kantorovich operators to a function}
\end{figure*}

\begin{figure*}[htb!]
\begin{center}
\includegraphics[height=6cm, width=10cm]{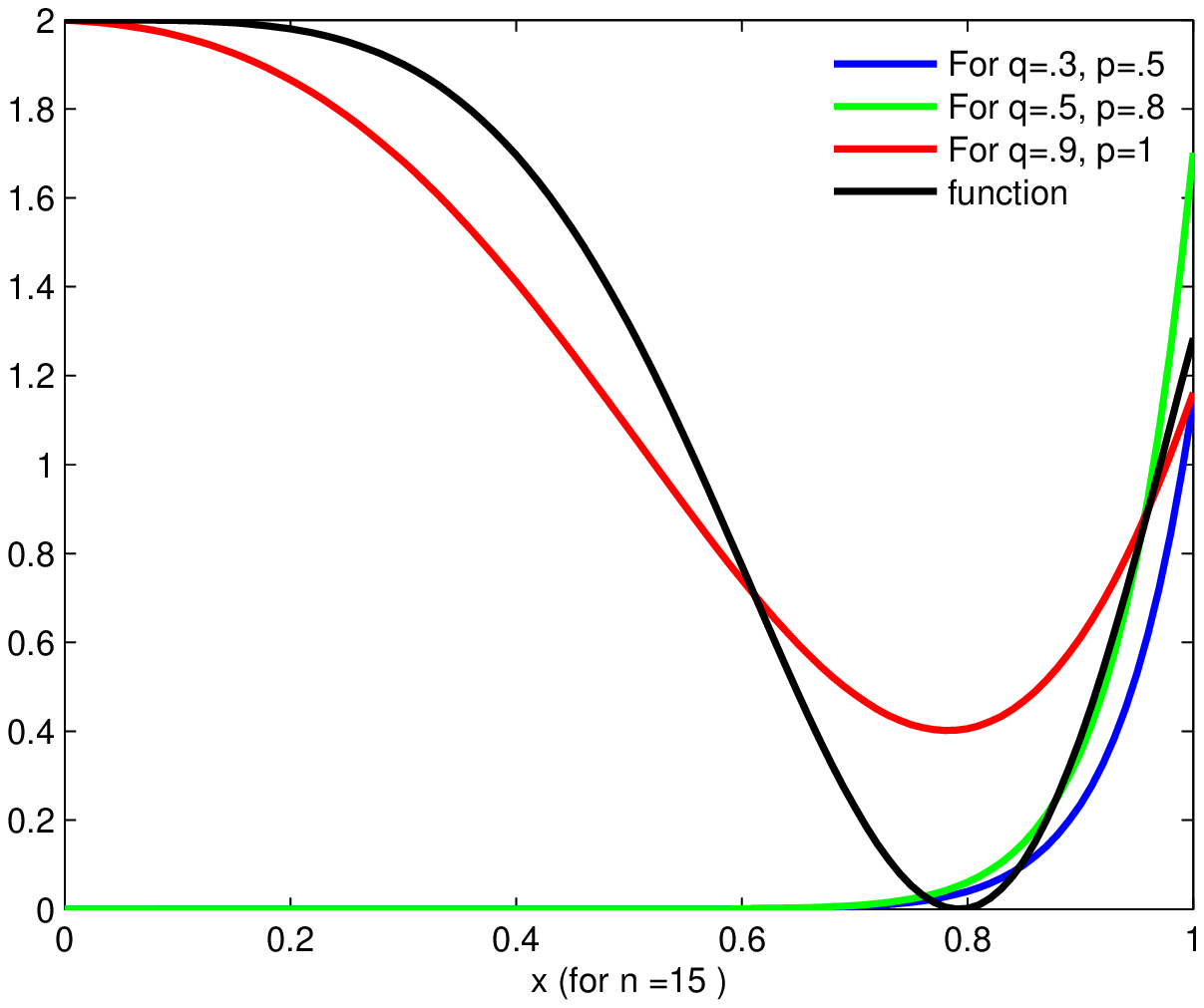}
\end{center}
\caption{Approximation by $(p,q)$-Bernstein-Schurer-Kantorovich operators to a function}
\end{figure*}

\begin{figure*}[htb!]
\begin{center}
\includegraphics[height=6cm, width=10cm]{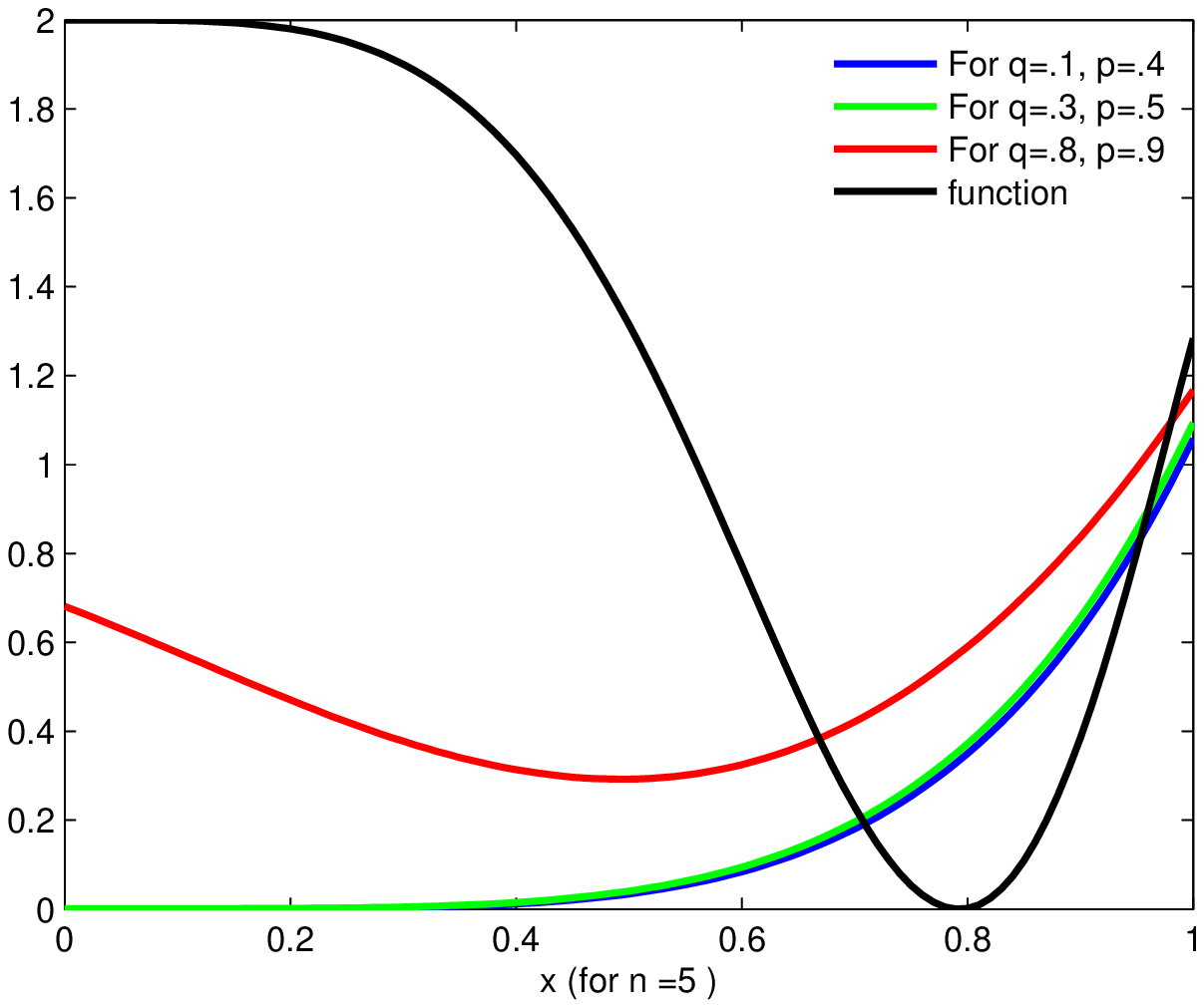}
\end{center}
\caption{Approximation by $(p,q) $-Bernstein-Schurer-Kantorovich operators to a function}
\end{figure*}

\begin{figure*}[htb!]
\begin{center}
\includegraphics[height=6cm, width=10cm]{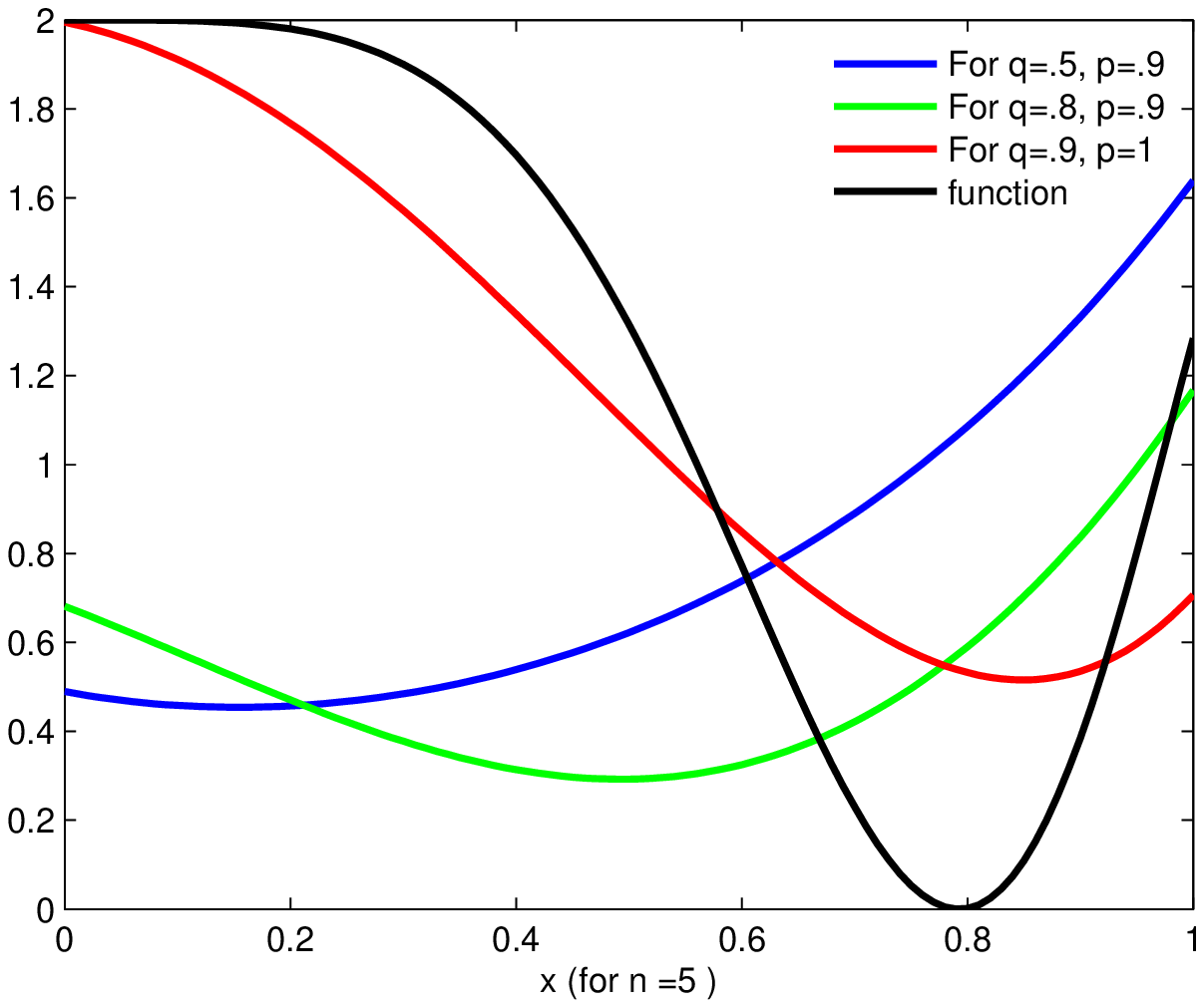}
\end{center}
\caption{Approximation by $(p,q) $-Bernstein-Schurer-Kantorovich operators to a function}
\end{figure*}


\newpage

\begin{thebibliography}{99}

\bibitem{ov} M.A. Ozarslan, Tuba Vedi, $q$-Bernstein-Schurer-Kantorovich Operators,
J. Inequal. Appl. , 256 (2013) 213-444.

\bibitem{gp1} A. Aral, V. Gupta and R.P. Agarwal, Applications of $q$-Calculus in Operator Theory,
Springer Science$\tiny{+}$Business Media New York, 2013.



\bibitem{alt} F. Altomare, M. Campiti, Korovkin type approximation theory
and its applications, de Gruyter Stud. Math. 17, Berlin, 1994.

\bibitem{ar1} A. Aral, O. Do\u{g}ru, Bleimann Butzer and Hahn operators
based on $q$-integers, J. Inequal. Appl., (2007) 1-12. Art. ID 79410.


\bibitem{brn} S.N. Bernstein, D\'{e}mostration du th\'{e}or\`{e}me de
Weierstrass fond\'{e}e sur le calcul de probabilit\'{e}s, Comm. Soc. Math.
Kharkow (2), 13 (1912-1913) 1-2.

\bibitem{bezier} P.E. Bezier, Numerical Control-Mathematics and applications, John Wiley and Sons, London, 1972.

\bibitem{chak} R. Chakrabarti and R. Jagannathan, A $(p,q)$-oscillator
realization of two parameter quantum algebras, J. Phys. A: Math. Gen., 24
(1991) 711-718.

\bibitem{dal} O. Dalmanoglu, Approximation by Kantorovich type q-Bernstein operators, in Proceedings
of the 12th WSEAS International Conference on Applied Mathematics, Cairo, Egypt (2007),
pp. 113-117.

\bibitem{dl} R.A. Devore, G.G. Lorentz, Constructive Approximation,
Springer, Berlin, 1993.

%

\bibitem{mah} M.N. Hounkonnou, J. D\'{e}sir\'{e}, B. Kyemba, $\mathcal{R}%
(p,q)$-calculus: differentiation and integration, SUT Journal of
Mathematics, Vol. 49, No. 2 (2013), 145-167.

\bibitem{kac} V. Kac, P. Cheung, Quantum Calculus, Springer-Verlag New York, 2002.

\bibitem{kan} L. V. Kantorovich, Sur certains d\'{e}veloppements suivant les
polyn\^{o}mes de la forme de S. Bernstein, I, II, C. R. Acad. URSS,
(1930), 563-568, 595-600.

\bibitem{pp} P. P. Korovkin, Linear operators and approximation theory,
Hindustan Publishing Corporation, Delhi, 1960.

\bibitem{jacob} J. Katriel, M. Kibler, Normal ordering for deformed boson
operators and operator-valued deformed Stirling numbers, J. Phys. A: Math.
Gen. 25 (1992) 2683-2691. Printed in the UK.

\bibitem{jag} R. Jagannathan, K. S. Rao, Two-parameter quantum algebras,
twin-basic numbers, and associated generalized hypergeometric series,
Proceedings of the International Conference on Number Theory and
Mathematical Physics, 20-21 December 2005.

\bibitem{lp} A. Lupa\c{s}, A $q$-analogue of the Bernstein operator, Seminar
on Numerical and Statistical Calculus, University of Cluj-Napoca, 9(1987)
85-92.


\bibitem{m3} N.I. Mahmudov and V. Gupta, On certain $q$-analogue of Sz\'{a}%
sz Kantorovich operators, J. Appl. Math. Comput., 37 (2011) 407-419.


\bibitem{ma1} M. Mursaleen, Asif Khan, Generalized $q$-Bernstein-Schurer
operators and some approximation theorems, Journal of Function Spaces and
Applications, Volume 2013, Article ID 719834, 7 pages
http://dx.doi.org/10.1155/2013/719834.

\bibitem{kl} Khalid Khan and D.K. Lobiyal, Bezier curves based on Lupas $(p,q)$-analogue of
Bernstein polynomials in CAGD, arXiv:1505.01810v1 [cs.GR].

\bibitem{mka1} M. Mursaleen, Khursheed J. Ansari, Asif Khan, On $(p,q)$-analogue of Bernstein Operators, (submitted).

\bibitem{mka2} M. Mursaleen, Khursheed J. Ansari, Asif Khan, Some Approximation Results by $(p,q)$-analogue of
Bernstein-Stancu Operators, (submitted).

\bibitem{mka3} M. Mursaleen, K. J. Ansari and Asif Khan, Approximation by a $(p,q)$-analogue of
Bernstein-Kantorovich Operators, arXiv: [1504.05887][math.CA]


\bibitem{pl} G.M. Phillips, Bernstein polynomials based on the $q$-integers,
The heritage of P.L.Chebyshev, Ann. Numer. Math., 4 (1997) 511-518.


 \bibitem{ph1} G. M. Phillips, A generalization of the Bernstein polynomials
based on the $q$-integers \emph{ANZIAMJ} 42(2000), 79-86.

\bibitem{hp} Halil Oruk, George M. Phillips,  q-Bernstein polynomials and Bezier curves, \emph{Journal of Computational and Applied Mathematics} 151 (2003) 1-12.

\bibitem{vivek} V. Sahai, S. Yadav, Representations of two parameter quantum
algebras and $p$,$q$-special functions, J. Math. Anal. Appl. 335 (2007)
268-279.

\bibitem{sad} P. N. Sadjang On the fundamental theorem of $(p,q)$-calculus
and some $(p,q)$-Taylor formulas, arXiv:1309.3934v1 [math.QA], 2013.


\end{thebibliography}
\end{document}